\documentclass[a4paper]{amsart}
\usepackage[T1]{fontenc}
\usepackage[latin1]{inputenc}
\usepackage{amssymb}
\usepackage{euscript}
\usepackage{amsxtra}
\usepackage[all]{xy}
\usepackage{enumerate}
\usepackage{mathrsfs}
\usepackage{amscd}
\usepackage{graphicx}
\usepackage[lite]{amsrefs}

\newcommand*{\Z}{{\mathbb Z}}      
\newcommand*{\C}{{\mathbb C}}      
\newcommand*{\Q}{\mathbb{Q}}

\newcommand*{\Comp}{{\mathbb K}}      
\newcommand*{\Bound}{{\mathbb B}} 

\newcommand*{\abs}[1]{\lvert#1\rvert}

\newcommand*{\cross}{\mathbin{\rtimes}}             
\newcommand*{\defeq}{\mathrel{:=}}

\newcommand*{\hot}{\mathbin{\hat{\otimes}}}

\newcommand*{\KK}{\mathrm{KK}}

\newcommand*{\K}{\mathrm{K}}

\newcommand*{\Eul}{\mathrm{Eul}}

\theoremstyle{plain}
\newtheorem{theorem}{Theorem}
\newtheorem{lemma}[theorem]{Lemma}

\theoremstyle{definition}
\newtheorem{definition}[theorem]{Definition}

\theoremstyle{remark}
\newtheorem{remark}[theorem]{Remark}

\newtheorem{example}[theorem]{Example}

\newcommand{\ctau}{C_\tau}

\newcommand{\sad}{s^{*}}

\newcommand{\Sigmaplus}{\Sigma^{+}}
\newcommand{\Pn}{\mathbf{P}}
\newcommand{\ind}{\mathrm{Index}}
\newcommand{\cstar}{C^{*}}
\newcommand{\Deltah}{\widehat{\Delta}}

\newcommand{\Psibar}{\dot{\Psi}}
\newcommand{\vs}{\mathbf{V}^{\mathrm{s}}}

\newcommand{\Hom}{\mathrm{Hom}}
\newcommand{\KKbullet}{\KK^{\bullet}}
\newcommand{\KKcbullet}{\KKbullet_\C}
\newcommand{\Kbullet}{\K_\bullet}
\newcommand{\Kcbullet}{\K_\bullet^{\C}}

\begin{document}

\title{Lefschetz numbers for $C^{*}$-algebras}

\author{Heath Emerson}
\email{hemerson@uvic.ca}

\address{University of Victoria\\ Victoria, BC}

\begin{abstract}

Using Poincar\'e duality, we formulate a formula of 
Lefschetz type which computes the Lefschetz number of 
an endomorphism of a separable, nuclear $\cstar$-algebra 
satisfying Poincar\'e duality and 
the Kunneth theorem.  (The Lefschetz number of an endomorphism 
is the graded trace of the induced map on $\K$-theory tensored 
with $\C$, as in the 
classical case.) We then examine endomorphisms of  
Cuntz-Krieger algebras $O_A$.  
An endomorphism has an invariant, which is a
permutation of an infinite set, and the contracting 
and expanding behavior 
of this permutation describes the Lefschetz number of 
the endomorphism. Using this description we derive a 
closed 
polynomial formula for the  Lefschetz number depending on 
the matrix $A$ and 
the presentation of the endomorphism.

\end{abstract}
 

\thanks{The author was supported by an NSERC Discovery grant}

\maketitle

\section{Introduction}
\label{sec:intro}

Suppose $A$ and $B$ are two separable, nuclear
 $\cstar$-algebras. To say that $A$ and $B$ 
are \emph{Poincar\'e dual} means that  
there is given a $\K$-homology class for $A\otimes B$ such 
that cup-cap product with this class induces an isomorphism 
between the $\K$-theory of $A$ and
 the $\K$-\emph{homology} of 
  $B$.   The homology class plays the role of the orientation class of 
  a compact manifold. 
The idea in this form is due to Alain Connes  (see \cite{Connes}). 
  Since the definition was invented,  
  quite a number of examples of Poincar\'e dual pairs have appeared 
  in the operator algebra literature, connected with dynamical systems, 
  foliations, hyperbolic groups, twisted $\K$-theory, $\cstar$-algebras 
of discrete groups with finite $B\Gamma$, \emph{etc}.

  The object of this note is 
  to propose a simple application of the existence of duality between a pair of
  algebras,  which runs roughly along the lines of a classical
    argument with de Rham cohomology and differential forms. 
  Suppose $\phi\colon X \to X$ is a smooth self-map of a compact, oriented
   manifold. Assume 
that $\phi$ is in general position with regard to fixed-points. Then $\phi$ induces 
a map on homology with rational coefficients, and its Lefschetz number is 
$$\mathrm{tr}_s (\phi_{*}) \defeq \mathrm{trace}(\phi_{*}\colon \mathrm{H}_{\mathrm{ev}} (X) \to 
\mathrm{H}_{\mathrm{ev}} (X)\bigr) - 
 \mathrm{trace}(\phi_{*}\colon \mathrm{H}_{\mathrm{odd}} (X) \to 
\mathrm{H}_{\mathrm{odd}} (X)\bigr).$$
The Lefschetz fixed-point theorem states that this number is equal 
to the number of fixed points of $\phi$ counted with appropriate multiplicities. The proof, which can be found in any textbook, involves ideas connected with Poincar\'e duality in de Rham theory: 
normal bundles, integration of forms, Thom classes, and so on. The Kunneth formula is a separate, additional ingredient. It is sometimes therefore said that the Lefschetz fixed-point formula \emph{follows} from Poincar\'e duality and the Kunneth formula.  

In this article, 
we are going to 
 formalize the exact connection between Poincar\'e duality and the Lefschetz 
fixed-point theorem in such a way as to apply to the category of $\cstar$-algebras, with 
$\K$-theory and $\K$-homology playing the role of ordinary homology and cohomology. 
Part of the proof of the classical Lefschetz theorem is absorbed into our statement, so 
that the classical theorem can be deduced from ours by a simple, essentially 
linear, index 
calculation.  We show via essentially straightforward formal calculations with 
  $\KK$-theory, that if one has a Poincar\'e duality with 
  `fundamental class' $\Delta \in \KK^n(A\otimes B, \C)$, and if one 
  has a morphism $f\in \KK(A,A),$ then the trace of the map  
  on $\K$-theory induced from $f$, 
  is equal to the result of a certain index pairing (see Section 2) 
  involving 
  $f$ and $\Delta$. This index pairing can in principle be computed in 
geometric terms, provided that the cycles underlying $f$ and $\Delta$ themselves  
 admit interesting, `geometric'  descriptions. Thus, to summarize, the Lefschetz number 
can be realized as a Kasparov product. Of course there is more than one such 
realization; in \cite{EmersonMeyer2} we pursue a similar idea to produce other 
kinds of identities in equivariant $\KK$-theory.

Of course the main merit of the observation is that one now has the possibility of 
proving analogues of the Lefschetz theorem in many different settings, provided one has 
available an interesting instance of noncommutative Poincar\'e duality. 

The significance of the classical Lefschetz formula tends to be explained in terms 
of the equality of a \emph{local} and a \emph{global} invariant. In connection with 
$\cstar$-algebras, this does not entirely make sense. What kind of Lefschetz formulas 
can we expect in connection with $\cstar$-algebras? One example, based on the 
abstract Lefschetz formula presented here, is worked out in \cite{EEK}. This 
involves proper actions of discrete groups on manifolds. The primitive ideal space 
of a cross product $C_0(X)\cross G$ in this situation is the extended quotient 
$$G\setminus \hat{X}, \; \text{where} \; \hat{X} \defeq \{(x,h) \in X\times G \; | \; h \in 
\mathrm{Stab}_G(x)\},$$
where $G$ acts on $\hat{X}$ by $g(x,h) = (gx, ghg^{-1})$, and 
which, as a set, identifies canonically with the primitive ideal space of $C_0(X)\cross G$ 
and inherits a corresponding hull-kernel topology. It is a bundle over the ordinary 
space $G\setminus X$ with fibre at $Gx$ the irreducible dual of $\mathrm{Stab}_G(x)$, 
but it is not Hausdorff. 
The Lefschetz formula for an 
automorphism of this situation has the corresponding shape: the geometric side of 
the formula involves fixed points in the ordinary space $G\setminus X$, and 
secondly, involves 
representation theoretic data for the isotropy of these fixed-points. 

The second purpose of this note is to consider the case of a pair of 
\emph{simple} algebras in duality, namely to pairs $A = O_A$, $B = O_{A^{T}}$ 
of Cuntz-Krieger algebras (see \cite{PK}).  Here, 
in contrast to the example of the previous 
paragraph, here there are no points at all. Given an endomorphism of 
  $O_A$ arising from certain geometric data, 
 we will solve the index problem on the geometric side of the formal 
  Lefschetz formula. 
 The endomorphisms with which we work correspond to $n$-tuples of 
continuous, partially defined homeomorphisms 
$$\varphi \colon Z \subset \Sigmaplus_A \to \Sigmaplus_A,$$
where $\Sigmaplus_A$ is the symbol space of sequences 
$(x_i)$ such that $A_{x_i, x_{i+1}} = 1$ for all $i$. The information 
involved in such an $n$-tuple can be summarized in a \emph{single} map  
on the countable set of paths in the graph corresponding to $A$. 
 The 
  geometric computation of the Lefschetz number turns out to depend, 
roughly, on the 
  difference between the number of strings whose length is shrunk by the 
  map, and the number of strings whose length is expanded by the map. 
This eventually leads to a description of the Lefschetz number roughly 
in the following terms: if we write 
$t_i = \sum s_\mu s_\nu^{*}$ for words $\mu, \nu$, 
for one of the images of the generators of $O_A$ under the endomorphism, 
then an appearance of 
$(\mu, \nu)$ with $\abs{\mu} \le \abs{\nu} $ contributes $+1$ to the Lefschetz 
number and $\abs{\mu} > \abs{\nu} +1$ contributes $-1$ (and there are 
no contributions when
 $\abs{\mu} = \abs{\nu} +1$). As a result of this decription we 
can, if we want, write down an explicit, closed formula for the Lefschetz 
number, which is a polynomial expression in the entries of the matrix 
$A$. 

This example thus shows that the `Lefschetz trick' results in interesting formulas even in what one might loosely refer to as a `very' noncommutative situation. 

The idea of formalizing the Lefschetz fixed-point theorem's proof using Poincar\'e 
duality (and the Kunneth theorem) in order to work in a more general context, 
is due to Andr\'e Weil, though 
not of course in connection with $\cstar$-algebras and $\K$-theory. It was used by him
 in connection with the so-called Weil conjectures 
 (see the Appendix to 
  \cite{Hartshorne}). So in this sense, we have rediscovered an old trick. However, 
even so it seems worth making it explicit in the operator algebraic context in view of 
the variety of 
Lefschetz-type formulas one can reasonably hope to achieve 
 by using $\cstar$-algebras and $\KK$-theory, which embrace such a 
wide selection of geometric situations.

  \section{The abstract Lefschetz theorem} 

Kasparov's $\KK$-theory is a realization of an additive, $\Z/2$-graded 
category with objects $\cstar$-algebras and morphisms $A \to B$ 
the elements of $\KKbullet (A,B)$, defined as a quotient of a certain 
set of cycles, by a certain equivalence relation (see \cite{Kasparov}.) 

In addition to its structure of an additive category, $\KK$ is a
 symmetric monoidal category with unit object the $\cstar$-algebra 
$\C$ and bifunctor given by the tensor product of $\cstar$-algebras on 
objects and the "external product" 
\begin{equation}
\KK^\bullet(D_1, D_1') \times \KK^\bullet (D_2, D_2') \to \KK^\bullet (D_1\otimes D_2, 
D_1'\otimes D_2'),\; \; \; (f_1, f_2) \mapsto f_1\hot_\C f_2. 
\end{equation}
on morphisms, and the flip 
$\Sigma \colon A\otimes B \to B\otimes A$ inducing the braiding. 


 The interaction between the flip, the monoidal structure, and the grading 
in $\KK$  is summarized by the following diagram, which 
\emph{graded commutes}, for all $D_i$, $D_i'$: 
$$
  \xymatrix{
     \KKbullet (D_1,  D_1') \times \KKbullet(D_2, D_2') \ar[d]_{\text{flip}}\ar[r]^{\;\;\;\;\;\;\;\hot_\C}  &
    \KKbullet (D_1\otimes D_2, D_1'\otimes D_2') \ar[d]^{\text{flip}} \\
    \KKbullet (D_2, D_2')\times \K_\bullet (D_1, D_1')\ar[r]^{\hot_\C} &
    \KKbullet (D_2\otimes D_1, D_2'\otimes D_1')}$$
In other words,  
\begin{equation}
\label{gradedcommute}
f_1\hot_\C f_2 = (-1)^{\partial f_1\partial f_2 } [\Sigma]\hot_{D_2\otimes D_1}
(f_2\hot_\C f_1) \hot_{D_2'\otimes D_1'} [\Sigma]
\end{equation}
for all $f_1\in \KKbullet (D_1, D_1')$ and $f_2 \in \KKbullet (D_2, D_2')$. 

Of course a category with similar properties is the category of complex 
$\Z/2$-graded vector spaces and vector space maps, where for
 the monoidal structure we use 
graded tensor product of vector spaces, and for the braiding we use
 the \emph{graded flip} 
 $$\Sigma^{\text{s}} (a\hot_\C b ) \defeq (-1)^{\partial a \partial b} b\hot_\C a$$ instead of the ordinary flip. The action of a linear transformation $T_1\hot_\C T_2 $ on $V_1\hot_\C V_2$, where $T_i\colon V_i \to V_i'$, $V_i$, $V_i'$ graded vector spaces, is defined by 
$$(T_1\hot_\C T_2) (a\hot_\C b) \defeq (-1)^{\partial T_1\partial b}T_1(a)\hot_\C T_2(b).$$
Then a short calculation shows that 
depending on the fact that $\partial x_1\partial x_2 + \partial T_1\partial x_2 + (\partial T_1 + \partial x_1) (\partial T_2 + \partial x_2)= 
\partial T_1\partial T_2 + \partial x_1\partial x_2$ mod $(2)$ shows that 
the monoidal structure on $\vs$ is also graded commutative, in the 
sense described above for $\KKbullet$.

These definitions ensure that 
the $\K$-theory functor 
$\KK \to \vs$, 
$$A \mapsto \Kcbullet (A), \;f\in \KKcbullet (A,B) \mapsto f_*\colon \Kcbullet (A) \to \Kcbullet (B),$$
which associates to a $\cstar$-algebra $A$ the complex, $\Z/2$-graded 
vector space $\Kcbullet (A) \defeq \K_\bullet (A)\otimes_{\Z}\C$, 
is compatible with the symmetric monoidal structures on each category, 
at least on a bootstrap category $\mathcal{N}$ (the Kunneth theorem) where 
 it is an isomorphism (the Universal Coefficient theorem.) 




In order to illustrate these facts in a concrete way we prove the following simple 
lemma. 

\begin{lemma}
\label{lem:basic} Suppose $c= \sum a_i\hot_\C b_i \in \Kbullet (A\otimes B)$ 
is written as a sum with 
$a_i, b_i$ homogeneous. \footnote{An element of a graded set 
is \emph{homogeneous} if it has a 
definite degree.} Let 
$f\in \KKbullet (A,A')$ and $g\in \KKbullet (B,B')$ be homogeneous.
 Then 
$$c\hot_{A\otimes B} (f\hot_\C g)  
= \sum (-1)^{\partial b_i\partial f} (a_i\hot_A f) \hot_\C (b_i\hot_B g) \in \Kbullet (A'\otimes B').$$
\end{lemma}

\begin{proof}
Suppressing subscripts, suppose $a\in \KKcbullet (\C, A)$, $b\in \KKcbullet (\C, B)$, 
$f$ and $g$ as above. Then  
\begin{multline}
(a\hot_\C b)\hot_{A\otimes B} (f\hot_\C g) \\
= a\hot_A (1_A\hot_\C b) \hot_{A\otimes B} (f\hot_\C 1_B) \hot_{A'\otimes B} (1_{A'}\hot_\C g)\\
= a\hot_A \Sigma_* (b\hot_\C 1_A) \hot_{A\otimes B} (f\hot_\C 1_B) \hot_{A'\otimes B} (1_{A'}\hot_\C g).
\end{multline}
Since $\Sigma_*(f\hot_\C 1_B) = \Sigma_*(1_B\hot_\C f)$, we can write the above 
\begin{equation}
\label{greenbananas}
= a\hot_A (b\hot_\C 1_A) \hot_{B\otimes A} (1_B\hot_\C f) \hot_{B\otimes A'} \Sigma^{*} (1_{A'}\hot_\C g).
\end{equation}
Using graded commutativity we have  
\begin{multline}
(b\hot_\C 1_A) \hot_{B\otimes A} (1_B\hot_\C f) = b\hot_\C f 
= (-1)^{\partial b \partial f} \Sigma_*(f\hot_\C b) \\ = (-1)^{\partial b \partial f} \Sigma_*\bigl( f\hot_{A'} (1_{A'}\hot_\C b)\bigr).
\end{multline}
Putting this into 
\eqref{greenbananas} and moving the flip across the tensor product gives 
\begin{multline}
= (-1)^{\partial b \partial f} a\hot_A f\hot_{A'} (1_{A'} \hot_\C b) \hot_{A'\otimes B} (1_{A'} \hot_\C g) \\ 
= (-1)^{\partial b \partial f} (a\hot_A f) \hot_\C (b\hot_B g),
\end{multline}
as required. 

\end{proof}

We next state the essential definition of this note (see \cite{Connes}, 
\cite{PK}, \cite{Emerson}, \cite{BMRS}.)

\begin{definition}
\label{pd}
\rm
Let $A$ and $B$ be $\cstar$-algebras. Then \emph{$A$ and $B$ are dual in $\KK$} (with a dimension shift of 
$n$) if there 
exists $\Delta \in\KK^{n}(A\otimes B, \C)$ such that the composition 
\begin{multline}
\label{pdmap}
\KKbullet (D_1, D_2\otimes A) \stackrel{- \hot_\C 1_B}{\longrightarrow} \KK^{\bullet}(D_1\otimes B, D_2\otimes A\otimes B) \\
\stackrel{\hot_{A\otimes B} \Delta}{\longrightarrow} \KK^{\bullet + n}(D_1\otimes B, D_2)
\end{multline}
is an isomorphism for every $D_1$, $D_2$. We call $\Delta$ the \emph{fundamental class} of the duality. 
\end{definition}

Suppose $A$ and $B$ are dual with class $\Delta$. In the above notation, set 
$D_1 = \C$ and $D_2 = B$. Then there is a unique class $\Deltah' \in \KK^{-n}(\C, B\otimes A)$ such that 
the isomorphism 
\begin{equation}
\label{yoneda}
\KK^{-n}(\C, B\otimes A) \stackrel{\cong}{\longrightarrow} \KK^{0}(B,B)
\end{equation} 
carries $\Deltah'$ to $1_B$. We call $\Deltah'$ the \emph{dual fundamental class}. By definition, we have 
the equation 
\begin{equation}
\label{firstequationofpd}
(\Deltah'\hot_\C 1_B) \hot_{B\otimes A\otimes B} (1_B\hot_\C \Delta) = 1_B.
\end{equation}
A simple computation shows that the map 
\begin{equation}
\label{secondpdmap}
\KKbullet (D_1\otimes B, D_2) \stackrel{\hot_\C 1_A}{\longrightarrow} \KKbullet (D_1\otimes B\otimes A, D_2\otimes A) 
\stackrel{\Deltah' \hot_{B\otimes A}}{\longrightarrow} \KK^{\bullet -n}(D_1, D_2\otimes A).
\end{equation}
is an inverse to \eqref{pdmap}. We obtain a second equation 
\begin{equation}
\label{secondequationofpd}
(1_A\hot_\C \Deltah' )\hot_{A\otimes B\otimes A} (\Delta \hot_\C 1_A) = 1_A .
\end{equation}

If one prefers to arrange things in a different logical pattern, one can start with a pair of 
classes $\Delta$ and $\Deltah'$ and insist that they satisfy the equations 
\eqref{firstequationofpd} and \eqref{secondequationofpd}. Then the map 
as in \eqref{pdmap}
can be shown to be an isomorphism with inverse 
\eqref{secondpdmap}.

\begin{remark}\label{conventions}
\rm
\label{conventions}
In the above notation, 
\begin{multline}
\label{easycomputation}
(1_A\hot_\C \Deltah') \hot_{A\otimes B\otimes A} (\Delta\hot_\C 1_A)\\
 = (-1)^{n}\, \bigl( 
\Sigma_*(\Deltah') \hot_\C 1_A \bigr)\hot_{A\otimes B\otimes A} \bigl( 1_A\hot_\C \Sigma^{*} (\Delta)\bigr).
\end{multline}
In \cite{Emerson} the definition of Poincar\'e duality involved classes 
$\Deltah \in \KK^{-n}(\C, A\otimes B)$ and $\Delta^n(A\otimes B, \C)$ satisfying 
appropriate equations. To connect our current discussion with that definition, 
set $\Deltah = \Sigma_*(\Deltah')$. Then by \eqref{easycomputation}, 
the analogues of equations \eqref{firstequationofpd} and \eqref{secondequationofpd}
are 

\begin{multline}
\label{otherequationsofpd}
(\Sigma_*(\Deltah)\hot_\C 1_B) \hot_{B\otimes A\otimes B} (1_B\hot_\C \Deltah) = 1_B,\\
(\Deltah\hot_\C 1_A) \hot_{A\otimes B\otimes A} \bigl(1_A\hot_\C \Sigma^{*}(\Deltah)\bigr) = (-1)^{n}\, 1_A
\end{multline}
which is as in \cite{Emerson}. 

Notice also that the roles of $A$ and $B$ are 
symmetric when $n$ is even and 
anti-symmetric when $n$ is odd. 
\end{remark}

Given $A$ and $B$ dual as above, define 
a $\Z$-bilinear map 
\begin{equation}
\label{eq:poincarepairing}
\K_\bullet (A) \times \K_\bullet (B) \to \Z, \;\; (x\, | \, y) \defeq   y \hot_B \hat{x},
\end{equation} 
where 
$\hat{x}$ denotes the Poincar\'e dual of $x$. 

\begin{lemma}
With the Poincar\'e duality pairing defined in \eqref{eq:poincarepairing},
$$ (x\, | \, y) = (-1)^{\partial x \partial y}(x\hot_\C y )\hot_{A\otimes B} \Delta $$
for homogeneous elements $x\in \K_\bullet (A)$, $y \in \K_\bullet (B)$.

\end{lemma}

\begin{proof}
Expanding the definitions, we have 
\begin{multline}
y\hot_B \hat{x}= y\hot_B  (x\hot_\C 1_B) \hot_{A\otimes B} \Delta = y\hot_B \Sigma_*(1_B\hot_\C x)\hot_{A\otimes B} \Delta \\ = 
y\hot_B (1_B\hot_\C x) \hot_{A\otimes B} \Sigma^{*}(\Delta) = (y\hot_\C x) \hot_{B\otimes A} \Sigma^{*}(\Delta) 
\\ = (-1)^{\partial x \partial y} \Sigma_*(x\hot_\C y) \hot_{B\otimes A} \Sigma^{*}(\Delta) = (-1)^{\partial x \partial y}(x\hot_\C y) \hot_{A\otimes B} \Delta.
\end{multline}
\end{proof}

Tensoring with the complex numbers we obtain a 
duality pairing  
$$(\; | \; )\colon \Kcbullet (A) \times \Kcbullet (B) \to \C.$$ This pairing 
is non-degenerate if $B$ satisfies the Universal Coefficient theorem. It 
is supported on $\{(x,y) \; | \; \partial (x) + \partial (y) = n\}$. 

Now note that if
 $A$ and $B$ are Poincar\'e dual, then $\Kbullet (A)$ and $\Kbullet (B)$ are 
finitely generated abelian groups (and for the same reason, if
 $A$ and $B$ are dual in $\KK_{\C}$ then 
$\Kcbullet (A)$ and $\Kcbullet (B)$ have finite rank.)

By elementary 
methods one can thus find a basis $(x_{\epsilon, i })$ for $\Kcbullet (A)$ and a dual basis 
$(x^{*}_{n-\epsilon, j}) $ for $\Kcbullet (B)$ with respect to $(\; | \; )$, \emph{i.e.} so 
that we have 
\begin{equation}
(x_{\epsilon , i}\, | \, x^{*}_{\eta , j}) = \delta_{\eta , n-\epsilon} \delta_{ij}.
\end{equation}

\begin{lemma}
In terms of the bases $(x_{\epsilon , i})$ and $(x_{\eta , j}^*)$, the class 
$\Deltah' $ is given by 
$$\Deltah' = \sum_{i,\epsilon} (-1)^{n- \epsilon} x_{n-\epsilon, i}^{*} \hot_\C x_{\epsilon , i}.$$
\end{lemma}

\begin{proof}
It suffices to show that the map 
\begin{equation}
\label{ourcurrentmap}
\KK_\C^{-n}(\C, B\otimes A) \stackrel{\hot_\C 1_B}{\longrightarrow}\KK_\C^{-n}(B, B\otimes A \otimes B) 
\\ \stackrel{\hot_{A \otimes B} \Delta}{\longrightarrow} \KK^{0}_\C(B,B)
\end{equation}
sends 
$\sum_{i,\epsilon} (-1)^{n- \epsilon} x_{n-\epsilon, i}^{*} \hot_\C x_{\epsilon , i}$ to the identity in $\KKcbullet (B,B)$. Since 
we are over $\C$, the UCT gives that 
$\KKcbullet (B,B) \cong \Hom_\C \bigl( \Kcbullet (B), \Kcbullet (B)\bigr)$. 
If $x\in \Kcbullet (B)$, and denoting our proposed formula for $\Deltah$ by $\widehat{\delta}$, 
then we have 
\begin{equation}
x\hot_B (\widehat{\delta} \hot 1_B) = (-1)^{\partial x \partial \widehat{\delta} } \widehat{\delta}\hot_\C x
\end{equation}
by Lemma \ref{lem:basic}. 
Hence the image of $\widehat{\delta}$ under \eqref{ourcurrentmap} sends 
$x\in \Kcbullet (B)$ to 
\begin{multline}
\label{wherexgoes}
(-1)^{n \partial x } \, \sum_{i,\epsilon} (-1)^{n- \epsilon} \, \bigl( x_{n-\epsilon, i}^{*} \hot_\C x_{\epsilon , i}\hot_\C x\bigr)  \, \hot_{B\otimes A\otimes B} (1_B\otimes \Delta)
\\ = (-1)^{n\partial x  }\, \sum (-1)^{n- \epsilon}\, x_{n-\epsilon, i}^{*}\, \cdot \, \bigl( ( x_{\epsilon ,i} \hot_\C x )\hot_{A\otimes B}  \Delta\bigr) \\= 
 \sum (-1)^{n \partial x + (n-\epsilon)  + \epsilon \partial x}x_{n-\epsilon, i}^{*}\, \cdot \,(x_{\epsilon ,i} \, | \, x).
\end{multline}
Now setting $x = x^{*}_{\gamma, j}$, each term vanishes save when 
$\epsilon = n-\gamma$, in which case the sign is $(-1)^{n\gamma + \gamma + (n-\gamma)\gamma} = +1$. 

\end{proof}

With these preliminaries out of the way, we can now state and prove the formal Lefschetz 
theorem for Poincar\'e dual pairs of $\cstar$-algebras alluded to in the introduction. 

Suppose we have a duality with fundamental classes $\Deltah' \in \KK^{-n}(B\otimes A, \C)$ and 
$\Delta \in \KK^{n}(A\otimes B, \C)$. Let 
$f\in \KK(B,B)$. Define 
$$\mathbf{Ind}(\Delta, f) \defeq  \Deltah' \hot_{B\otimes A} (f\hot_\C 1_A)
\hot_{B\otimes A} \Sigma^{*}\Delta \in 
\KK(\C, \C) \cong \mathbb{Z}.$$
As the notation suggests, this `index' only depends on $f$ and $\Delta$ subject to the 
condition that $\Delta$ implement a Poincar\'e duality. However, the way 
$\mathbf{Ind}$ is defined involves both $\Delta$ and the dual class 
$\Deltah'$, so that if one changes the duality, it is easy to check that 
two cancelling changes are 
introduced into $\mathbf{Ind}$, so that $\mathbf{Ind}(\Delta ,f)$ does not depend on 
the choice of $\Delta$. 

Motivated by the classical case, we define 
the \emph{Lefschetz number} $\mathrm{Lef}(f)$ 
of 
$f\in \KK^0(B,B)$ in the standard way as the graded trace 
\begin{multline}
\mathrm{Lef}(f) \defeq \mathrm{tr}_s (f_*\colon \Kcbullet (B) \to \Kcbullet (B))
\\ \defeq \mathrm{trace}_\C\bigl( f_*\colon \K_0^\C(B) \to \K_0^\C (B) \bigr) - 
\mathrm{trace}_\C\bigl( f_*\colon \K_1^\C (B) \to \K_1^\C (B).\bigr).
\end{multline}
of $f$ acting on the complexified $\K$-theory of $B$.  

\begin{theorem}
\label{thm:abstractlef}Let $A$ and $B$ be $\cstar$-algebras satisfying the Universal coefficient theorems and 
the Kunneth theorem. Suppose that $A$ and $B$ are dual with fundamental class 
$\Delta \in \KK^{n}(A\otimes B, \C)$ and dual class $\Deltah' \in \KK^{-n}(\C, B\otimes A)$. Then 
for any $f\in \KK^{0} (B,B)$, the Lefschetz number of $f$ is equal to the 
index   
\begin{equation}
\mathrm{Lef}(f) = \mathbf{Ind} ( \Delta , f).
\end{equation}
In particular, the zeta function
$$\zeta_f (t) \defeq \sum_{n=0}^{\infty} \mathbf{Ind} (\Delta, f^{n})t^{n}.$$
is rational. 
\end{theorem}

\begin{proof}
Let $f\in \KK^{0}(B,B)$. We can write 
$f_*(x^{*}_{\epsilon , i} ) = \sum_r f^{\epsilon}_{ir} x^{*}_{\epsilon , r}$. Hence 
\begin{equation}
(1_A\hot_\C f)_* (\Deltah) = \sum (-1)^{n-\epsilon }f_{ir}^{n-\epsilon} \, x^{*}_{n-\epsilon ,r}\hot_\C x_{\epsilon, i}.
\end{equation} 
Applying the flip gives 
\begin{equation}
\Sigma_*(1_A\hot_\C f)_* (\Deltah) 
= \sum (-1)^{n-\epsilon  + \epsilon(n-\epsilon)}\, 
f_{ir}^{n-\epsilon} \, x_{\epsilon, i} \hot_\C x^{*}_{n-\epsilon ,r}.
\end{equation} 
Finally, pairing this expression using $\Delta$ gives 
\begin{multline}
<\Sigma_*(1_A\hot_\C f)_* (\Deltah), \Delta> 
= \sum (-1)^{n-\epsilon  + \epsilon(n-\epsilon) }f_{ir}^{n-\epsilon} \,  ( x_{\epsilon, i}\hot_\C x^{*}_{n-\epsilon ,r}) \hot_{A\otimes B} \Delta \\
= \sum (-1)^{n-\epsilon  }f_{ir}^{n-\epsilon} \, ( x^{*}_{n-\epsilon ,r}\, | \,  x_{\epsilon, i}) \\
= \sum (-1)^{n-\epsilon}f_{ii}^{n-\epsilon} = \mathrm{tr}(f^{0}_*) - \mathrm{tr}(f^{1}_*) = \mathrm{tr}_s (f_*)
\end{multline} 
as required. 

The statement regarding the zeta function is an elementary consequence, 
see \cite{Hartshorne}. 

\end{proof}

Finally, we note that it is rather natural to call the Lefschetz number of the 
identity morphism $1_B \in \KK^{\bullet}(B,B)$ the \emph{Euler characteristic} of 
$B$; it is the difference in ranks of $\K_0(B)$ and $\K_1(B)$, and by our 
formal Lefschetz theorem it is the index 
\begin{equation}
\label{eq:GaussBonnet}
\Eul_B = <\Deltah , \Delta>,
\end{equation}
which is a sort of formal Gauss-Bonnet theorem. 

\begin{example}
\rm
Let $A $ be the $\cstar$-algebra of sections 
$\ctau (X)$ of the Clifford algebra of a compact manifold $X$, and $B=C(X)$. The 
best-known example of $\K$-theoretic Poincar\'e duality is in this situation. 
The class $\Delta$ is represented by 
the unbounded self-adjoint operator $D \defeq d+d^{*}$ acting on the bundle 
$\Lambda^{*}_\C (X)$ of differential forms on $X$, where $d$ is the de Rham differential,
 and the additional datum of the  
Clifford multiplication 
$$\ctau (X) \otimes C(X) \to \Bound( L^{2}(\Lambda^{*}_\C(X))\bigr).$$
The class $\Deltah$
involves a Clifford multiplication by an appropriate vector field on $X\times X$, acting 
on a submodule of $C(X)\otimes \ctau (X)$, but the important point  
is that this vector field vanishes on the diagonal. It is immediate that when we take the Kasparov 
product of $\Delta$ and $\Deltah$, we get simply the operator $D$ acting on 
forms $L^{2}(\Lambda^{*}_\C(X))$. Hence \eqref{eq:GaussBonnet} says 
that $$\Eul_X = \mathrm{Index} (D_{\mathrm{dR}}),$$
with $D_{\mathrm{dR}}$ the de Rham operator on $X$. See \cite{EmersonMeyer1} 
for a closely related computation. 

It is also a simple manner to deduce the classical Lefschetz fixed-point formula 
in the same way. The fact that the class $\Deltah$ is supported in a neighbourhood 
of the diagonal in $X\times X$ means that if we twist $\Deltah$ by a smooth map 
whose graph $X\to X\times X$ is transverse to the diagonal, and pair with $\Delta$, the result is supported 
in a small neighbourhood of the fixed-point set of the map. The latter is a discrete set. Thus, the formal 
Lefschetz theorem gets us this far, and to finish the computation we need to carry 
out a local index computation. (See \cite{EEK} for the details, in a more general 
context.) 
\end{example}

The reader wanting other simple examples may wish to consider an automorphism 
of a finite group. A pleasant noncommutative Lefschetz formula for this situation can be 
deduced using the formal Lefschetz formula. This formula  
gives the well-known relationship between the number of fixed points of the induced map 
$\hat{\zeta}\colon \widehat{G} \to \widehat{G}$ on the (finite) space of irreducible 
representations, and the number of `$\zeta$-twisted conjugacy classes' in $G$. The 
reference \cite{EEK} also contains this result.

In the remainder of this note, we are going 
to work out a highly noncommutative example (the algebras $A$ and $B$ 
are \emph{simple}). The merit of considering an example like this is that we do get a 
genuinely new equality of invariants -- a genuinely new Lefschetz theorem. The difficulty with 
the example is that it is not so easy to see what its meaning is. So that it is helpful 
to have the classical examples at hand, for comparison.

\section{Example -- endomorphisms of Cuntz-Krieger algebras}

Let $O_A$ (see \cite{CK1}, \cite{CK2}) be the Cuntz-Krieger algebra with (irreducible) 
matrix $A$, the universal $\cstar$-algebra generated by $n$ nonzero 
partial isometries 
$s_1, \ldots , s_n$  such that 
$$\sum A_{ij}s_js_j^{*} = \sad_is_i.$$
 We are going to illustrate the formal Lefshetz theorem given in the previous section by proving an 
analogue of the Lefschetz fixed-point theorem for (certain) endomorphisms of $O_A$.

We first remind the reader of the following theorem of Cuntz and 
Krieger. 

\begin{lemma}
\label{lem:ktheory} (see \cite{CK1}). The group $\K_0(O_A)$ is $\Z/(1-A^{t})\Z$, and the group $
\K_1(O_A)$ is the quotient of $\Z/(1-A^{t})\Z$ by its torsion subgroup.  
\end{lemma}

To compute the Lefschetz number of an endomorphism of $O_A$, we must therefore 
split off the free part of the $\K_0$ group and compute the images of a set of generators, and 
similarily, find free generators for $\K_1$ and compute their generators.

\begin{example}
\label{mainxzero}
A standing numerical example will be the case 
$A= \begin{pmatrix} 1 & 1 & 0 \\ 1 & 1 & 1 \\ 0 & 1 & 1\end{pmatrix}$. 
By Lemma \ref{lem:ktheory},
 $\K$-theory of $O_A$ is infinite cylic in each of dimensions $0$ and $1$, 
with free generator $[s_1s_1^{*}]$ the class of the projection $s_1s_1^{*}$ 
in even degree, and free 
generator the class $[s_1+s_3^{*}]$ of the unitary $s_1 + s_3^{*}$ 
in odd degree. Note that 
 $[s_2\sad_2] = 0$ and $[s_3\sad_3] = -[s_1\sad_1]$ in $\K_0$.

\end{example}

Any (unital) endomorphism 
$\alpha\colon O_A \to O_A$ maps each generator $s_i$ to a partial 
isometry $t_i \in O_A$ such that $t_1, \ldots , t_n$ satisfy the 
same relations. Conversely, by the universal property of $O_A$, any 
choice of $t_1, \ldots , t_n$  satisfying the Cuntz-Krieger 
relations gives rise to a `symmetry' of $O_A$, \emph{i.e.} an endomorphism. 
A well-known family of is the 
periodic $1$-parameter family given by the circle action 
$$s_i \mapsto zs_i, \; \text{where} \; z\in S^1, \; i = 1, 2, \ldots , n.$$
These endomorphisms are, however, obviously homotopic to the identity, whence 
we cannot expect very interesting Lefschetz numbers (they will all be zero, since
 the Euler characteristic 
of $O_A$ is zero.) Instead we are interested in more combinatorially defined 
endomorphisms. 

The following 
definition is vaguely analogous to the assumption that one has an infinitely differentiable 
map, in the setting of the classical Lefschetz theorem. Let $\Sigmaplus_A$ denote the 
Cantor set of sequences $(x_i)$ in the graph $\Lambda$ determined by 
the matrix $A$. 

\begin{definition}
\label{def:smooth}
Let $Z\subset \Sigmaplus_A$ be an open subset and 
$\varphi\colon Z \to \Sigmaplus_A$ be a 
continuous map with domain $Z$.
 Then $\varphi$ is \emph{smooth} if
\begin{itemize}
\item $Z$ is a cylinder set. 
\item There exists 
$k \in \mathbb{N}$ and a map 
$\psi'\colon \Pn_{\le k} \to \Pn$ such that 
$
\psi(x_1, x_2, x_3, \ldots ) \\ = (\psi' (x_1, \ldots , x_k), x_{k+1}, x_{k+2}, \ldots ), 
$ for all $x = (x_1, x_2, x_3, \ldots ) \in Z$,
\end{itemize}
where $\Pn \defeq \{ (x_1, \ldots , x_m ) \; | \; A_{x_i, x_{i+1}} = 1, m \ge 0 \}$ is the 
set of finite allowable strings in the alphabet determined by $A$, and 
$\Pn_{\le k}$ is the set of strings of length at most $k$.

\end{definition}

We allow the empty string $\emptyset$. With this convention, the left shift 
$$\sigma_A\colon Z \defeq \Sigmaplus_A \to \Sigmaplus_A$$ is 
smooth, since $\sigma_A (x_1, x_2, \ldots ) = (\sigma_A' (x_1), x_2, x_3, \ldots )$ 
where $\sigma_A' (x) = \emptyset$ for every string of length $1$.

 Our definition is actually 
closer to the idea of a quasi-conformal map. Note that 
$\Pn$ is the vertex set of 
the tree $\tilde{\Lambda}$ which is the universal cover of $\Lambda$. As such, it 
admits a canonical path metric: it is `Gromov hyperbolic' as a metric space, and 
so has a Gromov boundary. 

\begin{lemma}
A smooth map $\varphi \colon Z \to \Sigmaplus_A$ is the boundary 
value of a quasi-isometry $\varphi '\colon Z'  \to \tilde{\Lambda}$, where $Z'$ is a subset 
of $\Pn$.  
\end{lemma}

\begin{proof}
Suppose we are given 
a smooth map $\psi \colon Z \to \Sigmaplus_A$ in the above sense, with
 $k$ and $\psi'$ as in the 
definition. We can take the cylinder set $Z$ to be of the form 
$Z = \{ x\in \Sigmaplus_A \; | \; \pi_N (x) \in F\}$, where $\pi_N \colon \Sigmaplus \to \Pn_N$
is the projection, and $F$ is a finite subset of $\Pn_N$, and where 
$N$ is larger than $k$. Now whether an infinite string is or is not in the domain of 
$\varphi$ only depends on the first $N$ letters. Geometrically, 
the set of infinite strings $x$ with first $N$ letters belonging to a given fixed,
finite set of finite strings, is a clopen set of $\Sigmaplus_A$, and is the 
closure in the compactification of the tree, of the set of \emph{finite} strings 
of length at least $N$, and with first $N$ letters in the given set. Hence 
the set $Z$ is the boundary values of a subset $Z'\subset \Pn$, \emph{i.e.}
$Z = \overline{Z'} \cap \partial \tilde{\Lambda}$, 
where $Z'$ is 
the set of finite strings of length at least $N$ with first $N$ letters in $F$. 

Assuming now that 
the finite string $(x_1, \ldots , x_m)$ is in $Z'$, whence that 
any boundary point $x = (x_1, \ldots , x_m, x_{m+1}, \ldots )$  is in $Z$, we have  
$$\psi (x) = (\psi ' (x_1, \ldots , x_k), x_{k+1}, \ldots , x_{m}, \ldots )$$ which  
says that the last letter of $\varphi ' (x_1, \ldots , x_k)$ is allowed to be followed 
by $x_{k+1}$. Therefore we can define $\psi' (x_1, \ldots , x_m) \defeq 
 (\psi ' (x_1, \ldots , x_k), x_{k+1}, \ldots , x_{m})$. It is clear that $\psi' $ is a 
quasi-isometry of the tree. It therefore extends to the boundary $\partial \tilde{\Lambda} = \Lambda$ 
 and it is 
clear that its boundary values gives precisely $\psi \colon Z \to \Sigmaplus_A$. 

\end{proof}

A typical `geometric' endomorphism of $O_A$ will be specified by the following 
definition. 

\begin{definition}
A \emph{geometric endomorphism} of $O_A$, where $A$ is $n$-by-$n$,
 shall refer to the data of a 
partition $\Sigmaplus_A = Z_1 \cup \cdots \cup Z_n$ of the symbol space, 
and an $n$-tuple $\Psi = (\psi_1, \ldots , \psi_n)$ of smooth homeomorphisms 
$\psi_i\colon W_i \stackrel{\cong}{\to} Z_i \subset \Sigmaplus_A$,  such that $W_i = \bigcup_{A_{ij}=1}  Z_i$. 

\end{definition}

It is clear that 
such partially defined maps  
determine elements of 
$O_A$; we define 
\begin{equation}
\label{presentation}
t_i \defeq \sum_{\mu \in W_i' , \, \abs{\mu} = k} s_{\psi_i' (\mu)} s_{\mu}^{*},
\end{equation}
where the summation is over the words of length $k$ in $W_i' $, with 
$W_i' \subset \Pn$ with $\overline{W'} \cap \partial \tilde{\Lambda} =W_i$ as 
explained above, and where 
$\psi_i'$ are the extensions of the $\psi_i$ to $\Pn$, and, of course, where 
$k$ is sufficiently large. 

Then the range projection of 
$t_i$ identifies, in the obvious sense, with the image $Z_i$ 
of $\psi_i$, and the 
cokernel projection identifies with the domain of definition $W_i$ 
of $\psi_i$. Hence, due to condition 3), 
we get an endomorphism 
$\alpha_{\Psi}(s_i) \defeq t_i$ of $O_A$. 

\begin{remark} 
The identity endomorphism corresponds to the evident partition 
with $Z_i = \{ x \; | \; x \; \text{begins with $i$}\}$ and 
$\psi_i (x) = (i,x)$ for $x\in W_i  \defeq \cup_{A_{ij} = 1} Z_j$. 
\end{remark}

From now on we will abuse notation and
 denote by the same letter 
the partially defined maps $\psi_i\colon \Pn \to \Pn$, and the maps 
$\psi_i \colon \Sigmaplus_A \to \Sigmaplus_A$. Of course 
there is ambiguity in the choice of the lifts $\psi\colon \Pn \to \Pn$, but 
we fix choices once and for all. Similarly we will write $W_i$ instead of 
$W_i' $ and $Z_i$ instead of $Z_i'$.

\begin{example}
\label{mainexone}
A good example of a geometric endomorphism 
for $A= \begin{pmatrix} 1 & 1 & 0 \\ 1 & 1 & 1 \\ 0 & 1 & 1\end{pmatrix}$ is  
\begin{multline}
\label{mainexampleone}
t_1 = s_1^{2}\sad_1\sad_2 + s_1s_2(\sad_2 )^{2} + s_2s_3^{2}\sad_3\sad_2 + s_2s_3s_2\sad_2\sad_3 + s_2\sad_1,  
\\ t_2 = s_3s_2$, $t_3 = s_3^{2}\sad_3.\end{multline}

The corresponding partition and triple of smooth maps is as follows.
 
\begin{itemize}
\item[1)] $Z_1$ is all sequences $(x_n)$ beginning with $1$ or $2$. 
\item[2)] $Z_2$ is all sequences $(x_n)$ beginning with $32$. 
\item[3)] $Z_3$ is all sequences $(x_n)$ beginning with $33$. 
\end{itemize}

\begin{itemize}
\item[\underline{$\psi_1$}] We require
 $\psi_1\colon W_1\defeq Z_1 \cup Z_3 \stackrel{\cong}{\longrightarrow} Z_1$. If a sequence 
begins with $1$, then we replace the initial $1$ by a $2$. If a sequence begins with 
$2$ then we replace the initial $2$ by a $1$, unless the second coordinate is $3$. 
In that case, we replace the initial $23$ by $233$. Finally, on sequences beginning with 
$32$, we replace the initial $32$ by $232$. Observe that the image of $\psi_1$ is 
all strings beginning with $2$ or $1$. 
\item[\underline{$\psi_2$}] We require $\psi_2\colon W_2 \defeq 
Z_1 \cup Z_2 \cup Z_3 = \Sigmaplus_A \stackrel{\cong}{\longrightarrow}
 Z_2$. We add $32$ to 
the beginning of any sequence. 
\item[\underline{$\psi_3$}] We require $\psi_3\colon W_3\defeq Z_2\cup Z_3 
\stackrel{\cong}{\longrightarrow} Z_3$.
 To any sequence beginning with $3$ we add an
additional $3$. 
\end{itemize}

\begin{remark} 
\label{rem:K_theory_computation}
The endomorphism $\alpha_\Psi$ above
 sends the range projection of $s_1$ to the range 
projection of $t_1$, which is $s_1\sad_1 + s_2\sad_2$. Hence
$(\alpha_\Psi)_*([s_1\sad_1]) = [s_1\sad_1] + [s_2\sad_2] = [s_1\sad_1]$, 
so the induced map $(\alpha_\Psi)_*\colon \K_0(O_A) \to \K_0(O_A)$ is 
the identity. To see the action 
on $\K_1 (O_A)$, one can check that the map 
$(t_1 + t_3^{*}  \; | \;\cdot \;) \colon \K_0(O_A) \to \Z$ induced from the 
Poincar\'e duality 
pairing (see the end of this section) vanishes identically. 
Therefore $[t_1 + t_3^{*}] = 0 \in \K_1(O_A)$ and 
$(\alpha_\Psi)_*\colon \K_1(O_A) \to \K_1(O_A)$ is the 
zero map.  So the 
Lefschetz number of $\alpha_\Psi \colon O_A \to O_A$ is equal to $1$. 
(In particular, $\alpha_\Psi$ is not an \emph{automorphism}.) 
\end{remark}

\end{example}

We now describe an invariant of any geometric endomorphism, which will 
be a 
\emph{single} partially defined map $\Pn \to \Pn$.

\begin{definition}
Let $\Psi = (\psi_1, \ldots , \psi_n)$ be a geometric endomorphism. Extend the 
$\psi_i$ to partially defined self-maps of $\Pn$. We let 
$\Psibar \colon \Pn \to \Pn$ be the partially defined map 
defined by 
$$\Psibar(x_1, \ldots , x_n)\defeq \psi_{x_n} (x_1, \ldots , x_{n-1})$$
if $(x_1, \ldots , x_{n-1}) \in \mathrm{Dom} (\psi_{x_n})$. 

\end{definition}

\begin{example} 
\label{mainextwo}
The map $\Psibar$ of Example \ref{mainexone} 
is define on paths of length $2$ by  
\begin{multline}
\Psibar (11) = 2,\;  \Psibar(21) = 1, \; \Psibar (12) = (321),\; \Psibar (22) = (322),\\  \Psibar (32) = (323),\; 
\Psibar (33) = (33),
\end{multline}

On paths of length $3$. 
\begin{multline}
\Psibar(111) = ( 21),  \Psibar( 121) = (22), \; \Psibar(211) = ( 11), \; 
\Psibar(221) = ( 12).
\end{multline}
\begin{equation}\Psibar (\star \star 2) = (32\star \star) \; \text{for any $(\star \star)$ allowable,}
\end{equation} 
\begin{equation}
\Psibar(323) =  (33),  \; \Psibar(333) = (333).
\end{equation}
Finally, on words of length $4$, $\Psi$ is defined by 
\begin{multline}
\Psibar(1111) = ( 211), \;
\Psibar( 1121) = (212),\;
\Psibar(1221) = (222),\; \\
\Psibar(2111) = ( 111), \; 
\Psibar(1211) = ( 221), \; 
\Psibar(2211) = ( 121), \\
\Psibar(3211) = ( 2321), \;
\Psibar(1121) = ( 212), \;
\Psibar(2121) = ( 112), 
\end{multline}

\begin{equation}
\Psibar( \star \star \star 2)= (32\star \star \star ), \; \text{for any $\star \star \star$ allowable,}
\end{equation}
and finally, 
\begin{multline}
\Psibar(3223) = (33223),\; 
\Psibar(3233) = (33233),\; 
\Psibar(3323) = (33323) , \; 
\\ \Psibar(3333) = (33333) , \;
\end{multline}

In these formulas, any string not mentioned is not in the domain. 

\end{example}

Say that two partially defined maps 
$\Psibar$ and $\Psibar$ are equivalent if $\Psibar = \Psibar'$ on sufficiently long strings. 
Then it is only the class of $\Psibar$ modulo $\sim$ which will matter to us for 
what is coming. Denote by 
 $[\Psibar]$ the class of $\Psibar$. 
We are going to associate to a geometric endomorphism
 an integer 
invariant. This invariant will only depend on the equivalence class 
$[\Psibar]$ and not on $\Psibar$ itself.  

By \emph{formal series} $\sum_{k=1}^\infty a_k$, where 
$a_k$ are real numbers, we will refer to the 
sequence of of its terms, modulo the equivalence relation 
$\sum_{k = 1}^\infty a_k \sim \sum_{k =1}^\infty b_k$ if $\sum_{k=1}^m a_k = \sum_{k =1}^m b_k $ for 
$m$ sufficiently large.  For example $1 + 2 + 3 + 4 + \cdots \sim 2 + 0 + 3 + 4 + \cdots$. 
The condition $\sim$ implies, obviously, that $a_k = b_k$ for large enough $k$, but it is stronger.

\begin{definition} 
\label{theindex}
\rm 
Let $\Xi \colon \Pn \to \Pn$ be
a partially defined bijection with finite propagation: that is, 
there exists $N \defeq 
\mathrm{Prop}(\Xi)$ such that 
$$\Xi (\Pn_k) \subset \bigcup_{\abs{l-k}\le N} \Pn_l, \; \; \text{for all $k$}.$$
Let
$\mathrm{Dom}(\Xi) \subset \Pn$ be its domain and 
$\mathrm{Im}(\Xi)$ its range. We set 
\begin{equation}
\mathrm{Index}_k (\Xi)  \defeq \mathrm{card}\, \bigl( \Pn_k \cap \mathrm{Im}(\Xi) \bigr)- 
\bigl( \Pn_k \cap \mathrm{Dom}(\Xi) \bigr).
\end{equation}
We let $\ind (\Xi)$ be the formal series   
\begin{equation}
\label{index}
\mathrm{Index} (\Xi) = \sum_{i = 1}^{\infty} \mathrm{Index}_i (\Xi).
\end{equation}
We show below that the index only depends 
on $[\Xi]$ and converges if $\Xi = \Psibar$ for some geometric endomorphism $\Psi $. 
\end{definition}

\begin{lemma}
\label{invariance}
The index has the following properties. 
\begin{itemize}
\item[1)]  If $\Xi$ and $\Xi '$ are 
two partially defined 
maps which agree on $\Pn_k$ for all sufficiently large $k$, then 
$\mathrm{Index}(\Xi) = 
\mathrm{Index}(\Xi')$ as formal series.   
Hence, $\ind$ is compatible with $\sim$. 
\item[2)] For any geometric endomorphism $\Psi $, 
$\ind_k(\Psibar)= 0$ for $k $ sufficiently large.  
Hence the formal series in \eqref{index} converges in this case.  
\end{itemize}

\end{lemma}

\begin{proof} 
We prove 2) first. For simplicity, we assume that a 
given partially defined map 
$\Xi$ has propagation at most $1$. We start by assuming that $\Xi$ has no strings in its domain of length 
$\le m-1$, for some $m \ge 2$. Now we remove any strings of length $m$ from the 
domain of $\Xi$. Let the new partially defined map be called $\Xi'$. We claim that 
the index (or more precisely the formal sum \eqref{index}) has not changed.  The   
index in dimension $m-1$ has clearly been shrunk by the number of elements in dimension $m$ 
which previously mapped to dimension $m-1$. Call this $a(m,m-1)$. Thus, 
$$\mathrm{Index}_{m-1}(\Xi') = \ind_{m-1}(\Xi) - a(m,m-1).$$ On the other hand, 
the domain in dimension $m$ has been reduced by $\mathrm{card}( \mathrm{Dom}(\Xi)\cap 
\Pn_m)$, while the image in dimension $m$ has been reduced by 
$a(m,m)$. So 
$$\mathrm{Index}_m(\Xi') =\ind_m(\Xi) -a(m,m) + 
\mathrm{Dom}(\Xi)\cap 
\mathrm{card}(\Pn_m) .$$
Finally, the image in dimension $m+1$ is reduced by $a(m,m+1)$. Meanwhile, the 
index in dimension $<m-1$ has not changed, nor has the index in dimensions $>m+1$, since 
$\Xi$ changes lengths of strings by at most $1$. So 
\begin{multline}
\ind (\Xi')  \\ = \ind(\Xi) -a(m,m-1)  -a(m,m) + -a(m,m) +  
\mathrm{card}(\mathrm{Dom}(\Psibar) \cap \Pn_m)\\  = \ind(\Xi).
\end{multline}
This proves the result. 

Now this means that for any $\Xi$ can have its domain 
successively shrunk by eliminating strings of length $1$, then $2$, and so on,  
without altering its index. The first assertion is now immediate, since without 
changing the index, we can alter both maps until they agree as partially defined maps. 
  
The second assertion is left to the reader (it follows from the analytic considerations discussed in the last section of the paper.) 
\end{proof}

\begin{example}
\label{classicfour}
Consider Example \ref{mainexone}, \ref{mainextwo}. 
The domain in dimension 
$1$ has $0$ elements in it. The image has $2$ elements in it. So $\ind_1(\Psibar) = 2$. 
The domain in dimension 
$2$ has 
$6$ elements in it and the image has $6$ elements in it. Hence $\ind_2(\Psibar) = 0$.
In dimension $3$ there are $13$ elements in the domain and $12$ in the image. So 
$\ind_3(\Psibar) = -1$. One checks that $\ind_k(\Psibar) = 0$ for $k >3$. Hence 
$$\ind(\Psibar) =  2 + 0 -1 = 1.$$

\end{example}

\begin{remark}\rm
Altering the domain of a $\Psibar$ on a finite piece is analogous to altering a map $f\colon M \to M$ up to 
homotopy, whilst retaining transversality.  The net effect on the fixed points (with signs) 
is zero. 
\end{remark}

\begin{remark}
 1) The identity morphism of $O_A$ corresponds to the partially defined map 
$\Psibar_{\mathrm{id}}$ 
with domain of definition the set of paths $(x_1, \ldots , x_n)$ such that 
$A_{x_n, x_1} = 1$, \emph{i.e.} the set of loops in the graph. The action of $\Psibar_{\mathrm{id}}$ is 
by shifting the parameterization of loops. In particular, $\mathrm{Index}(\Psi_{\mathrm{id}}) = 0$. 
2) If the graph corresponding to $A$ 
is \emph{complete}, then $\mathrm{Index}(\Psibar) = 0$ for \emph{every} $\Psibar$. 
This  
follows from the Lefschetz theorem.  
\end{remark}

The point about the index is that there is a lot of cancellation in the expression 
\eqref{index}. Taking into account this cancellation, we get a much more computable 
description of the index 

\begin{lemma}\rm 
\label{lem:main_technical}
Let $\Xi\colon \Pn \to \Pn$ be a partially defined homeomorphism with finite propagation.  
Let $m>0$. Define 
$$\gamma_m (\Xi) \defeq \sharp  \{ x\in \Pn\; | \; \abs{x} > m,
\; \abs{\Xi (x)} \le m \} \, -\, 
\sharp  
\{ x\in \Pn\; | \; \abs{x} \le m , \; \abs{\Xi (x)} >m \}.$$
Then 
$$\ind_1(\Xi) + \ind_2 (\Xi) + \cdots + \ind_m (\Xi) = \gamma_{m}(\Xi).$$
In particular,  if $\Xi = \Psibar$ for some geometric endomorphism $\Psi$,  then
$\gamma_m = \gamma_{m+1} =\cdots = \ind (\Psibar)$ for $m$ sufficiently large. 

\end{lemma}



\begin{proof}
Let $a(i,j)$ denote the number of strings of length $i$ which are mapped by $\Xi$ to 
strings of length $j$. Let $\delta(i,j) \defeq a(i,j) - a(j,i)$. Assume that $\Xi$ alters 
lengths of strings by at most $N$. Choose $k > 0$. 
By definition,
$$\ind_m(\Xi) = \sharp \, ( \mathrm{Im}(\Xi)\cap \Pn_m) - 
\sharp \,(\mathrm{Dom}(\Xi)\cap \Pn_m) .$$On the other hand, 
$\sharp \, \mathrm{Im}(\Xi)\cap \Pn_m) = \sum_{k = -N}^{N} a(m+k, m)$ while 
$\sharp \, (\mathrm{Dom}(\Xi)\cap \Pn) = \sum_{k =-N}^{N} 
a(m,m+k),$ whence 
$$\ind_m(\Xi) = \sum_{k = -N}^{N} \delta(m+k, m).$$
Of course $\delta(i,j) = -\delta (j,i)$. Hence when we take the (formal) sum 
\begin{equation}
\label{sum}
\ind (\Xi)  = \sum_{m=1}^{\infty} \sum_{k=-N}^{N} \delta(m+k, m),
\end{equation}
a term $\delta(i,j)$ appears exactly twice with opposite signs, if $i$ and $j$ are small enough relative to $m$. It follows that $$\sum_{k =1}^m \ind_k (\Xi) = \gamma_m (\Xi)$$
 because of telescoping. 
 
 The last assertion follows from Lemma \ref{invariance}.




\end{proof}

\begin{example}
For instance, in Example \ref{mainexone}, \ref{mainextwo},  
$m=3 $ 
is large enough, note $\mathrm{Prop}(\Psibar) \le 1$.   
There are $8$ strings of length $4$ mapped to strings of length $3$ and 
$7$ strings of length $3$ mapped to strings of length $4$, so 
 $$\ind(\Psibar) =  8-7 = 1.$$ 
\end{example}

Based on Lemma \ref{lem:main_technical}, we
 can give a poynomial formula for the index as follows. Fix $m$ large. 
Fix $j$. We count the number of strings of length $m+j$ (for $j = 1, 2, \ldots , N$, which 
are mapped to strings of length $\le m$. We refer to the presentation \eqref{presentation}. 
Fix $i$ and $\mu$ with $\abs{\mu} = k$. 
Suppose that $\abs{\psi_i (\mu)} \le \abs{\mu} - j +1$. 
Consider a string $w = (\mu, u)$ of length $m + j$, where 
$\abs{u} = m + j -\abs{\mu}$ is a string (path) from 
the terminus 
$t(\mu)$ of $\mu$ to $i$. Then this is mapped under $\Psibar$ to a string of 
length $\le m+j -1 -j+1 = m$. Hence for each such $i,\mu$ and $u$ we get a 
positive contribution to the index. For fixed $\mu$ and $i$ the number of 
possible $u$'s is equal to the number of paths of length $m+j-k$ 
from $t(\mu)$ to $i$, which equals 
 $A^{m+j-k}_{t(\mu) i}$. Hence the total positive contribution to the index is 
$$\sum_{i = 1}^{n}\sum_{j = 1}^{N}\sum_{\mu \in W_i, \, \abs{\psi_i (\mu)}
 \le \abs{\mu} - j + 1} A^{m+j-k}_{t(\mu) , i}.$$

For the negative contributions, for $j =0,1, \ldots , N-1$ fix $i$ and $\mu$ such 
that $\abs{\psi_i (\mu)} \ge \abs{\mu} + j + 2$. Then 
for each $w = (\mu, u)$ of length $m-j$, so that $u$ is a string from 
$t(\mu)$ to $i$ of length $m-j-k$, the length of $\Psibar (w)$ is 
$\ge m-j -1 + j+2 = m+1$. Hence we get a negative contribution to 
the index. Therefore the total negative contributions is  
$$\sum_{i = 1}^{n}\sum_{j = 0}^{N-1} \sum_{\mu \in W_i \, \abs{\psi_i (\mu)}\ge \abs{\mu} +j +2} A^{m-j-k}_{t(\mu),i}.$$
Therefore we get the following curious, completely explicit, polynomial 
formula for the index: it
 is given explicitly by the formula 
\begin{equation}
\label{eq:polynomial_formula}
\sum_{i = 1}^{n}\sum_{j = 1}^{N}\sum_{\mu \in W_i, \, \abs{\psi_i (\mu)}
 \le \abs{\mu} - j + 1} A^{m+j-k}_{t(\mu) , i} \, - \, 
\sum_{i = 1}^{n}\sum_{j = 0}^{N-1} \sum_{\mu \in W_i \, \abs{\psi_i (\mu)}\ge \abs{\mu} +j +2} A^{m-j-k}_{t(\mu),i},
\end{equation}
for any $N > \mathrm{Prop}(\Psibar) = \mathrm{max}_{i,\mu} \bigl( \abs{\mu} - \abs{\psi_i (\mu)} + 1\bigr)$ 
for any $m $ large enough. 

\begin{example}
For instance, in our main example the above formula with $k = 2$, 
$m = 3$, $N= 1$ gives 
\begin{multline}
\label{eq:example_index}
 \mathrm{Index}(\Psibar)  = (A^2_{11} + A^2_{21} + 
 A^2_{11} + A^2_{21} ) \\ - (A_{12} + A_{22} + A_{12} + A_{22} + A_{32} + A_{22} + A_{32})  = 8-7 = 1.
 \end{multline}
 \end{example}

We can now state our Lefschetz formula for Cuntz-Krieger algebra endomorphisms, at 
least those coming from simple combinatorics of generators and relations. 

\begin{theorem}(Lefschetz Theorem for Cuntz-Krieger endomorphisms.)
\label{mainlef} Let $\Psi \in G_A$ and $\alpha_\Psi \colon O_A \to O_A$ be the 
corresponding endomorphism. For sufficiently large \(m\) and \(N\), 
the 
Lefschetz number of $\alpha$ equals the index of $\Psi$. Thus, the 
Lefschetz formula 
\begin{multline}
\mathrm{Trace}_s \bigl((\alpha_\Psi)_* \colon \K_*(O_A)_\Q\to \K_*(O_A)_\Q\bigr) \\ = 
\sum_{i = 1}^{n}\sum_{j = 1}^{N}\sum_{\mu \in W_i, \, \abs{\psi_i (\mu)}
 \le \abs{\mu} - j + 1} A^{m+j-k}_{t(\mu) , i} \, - \, 
\sum_{i = 1}^{n}\sum_{j = 0}^{N-1} \sum_{\mu \in W_i \, \abs{\psi_i (\mu)}\ge \abs{\mu} +j +2} A^{m-j-k}_{t(\mu),i},
\end{multline}
holds. 
\end{theorem}

\begin{remark}
\label{eq:example_index}
The reader will clearly see the difference between computing 
the trace of an endomorphism using the polynomial formula 
\eqref{eq:polynomial_formula}, 
and by way of computing the action on \(\K\)-theory. Compare for instance the 
expression \eqref{eq:example_index} with the calculations 
in Remark \ref{rem:K_theory_computation}, which of course need the work of 
Cuntz and Krieger just to get off the ground. 
\end{remark}

To prove this, we need to show 
that $\mathrm{Index}(\Psibar) = \mathbf{Ind}(\Delta, [\alpha_\Psi])$ 
for an appropriate $\Delta$ inducing a duality. Kaminker and 
Putnam proved such a duality in 
\cite{PK}. We refer the reader to their paper for further details, and 
merely sketch the  computation here. Let $s_1, \ldots , s_n$ denote 
the generators for $O_A$ and $t_1, \ldots , t_n$ the generators for 
$O_{A^t}$. 

Define 
$H_A \defeq \ell^2(\Pn)$, where $\Pn$ is the set of strings, as above. 
Let 
$$S_i\colon H_A \to H_A, \; S_i(e_{w} ) \defeq A_{i\mathrm{o}(w)} e_{iw}, \;\;
R_j (e_w) \defeq A_{\mathrm{t}(w)j}e_{wj}.$$
Clearly $[S_i,R_j] = 0$, while $S_i^*, R_j] = 0$ modulo finite-rank operators. 
It is also easy to check that $\sum_j A_{ij} S_jS_j^* = S_i^*S_i$ modulo finite
rank operators, and similarly the $R_j$ satisfy the relations for 
$O_{A^t}$. 
We obtain the Busby invariant 
$$O_A \otimes O_{A^T} \to \Bound(H_A)/\Comp(H_A)$$
of an extension of $O_A\otimes O_{A^t}$ by the compact operators 
and hence (since $O_A$ is nuclear) a class in $\KK^1(O_A\otimes O_{A^t}, \C)$. 
Kaminker and Putnam prove that $\Delta$ induces a duality 
with dual class the element 
$w = \sum s_i\otimes t_i^*$. Then $ww^* = w^*w $ and each are projections. 
Therefore $w+ 1-ww^*$ is a unitary in $O_A\otimes O_A^t$ and so defines 
an $\Deltah$ of $\KK^1(\C, O_A, \otimes O_{A^t})$. Now suppose we have 
an endomorphism 
$$s_i \mapsto t_i \defeq \sum_{\mu \in W_i, \; \abs{\mu} = k} s_{\psi_i (\mu)}s_\mu^*.$$
Then under the endomorphism, $\Deltah$ is mapped to 
$$\sum_{i, \, \mu \in W_i, \; \abs{\mu} = k} s_{\psi_i (\mu)}s_\mu^*\otimes t_i.$$
To compute the pairing 
$$\mathbf{Ind} (\Delta, [\alpha]) =  <(\alpha_\Psi \otimes 1_{O_{A^t}})_*(\Deltah) , \Delta>$$
we need to compute the index of the obvious lift 
Fredholm 
index of $W_\Psi + (1-W_\Psi W_\Psi ^{*})$. However, it is clear that 
$W_\Psi W_\Psi^* $ is \emph{equal} to 
$W_\Psi^*W_\Psi$ on $\ell^2(\Pn_m)$ for $m\ge \mathrm{dmin}(\Psi) + 1$. 
Hence 
the sum 
$$\sum_{j=1}^{\infty} \mathrm{dim}\, \mathrm{ker}((W_\Psi)_{|_{\ell^2(\Pn_j)}}) 
- \mathrm{dim}\, \mathrm{ran}(W_\Psi) \cap \ell^2(\Pn_j))$$
converges, and evidently converges to the analytic index. 
Now we can regard $W_\Psi$ as the operator induced by the 
partial permutation $\Psibar$ of $\Pn$, in which point masses $e_w$ 
in the kernel of $W_\Psi$ correspond to words $w$ not in the 
domain of $W_\Psi$. We are now in the setting of our earlier 
discussion of partially defined maps $\Pn \to \Pn$, and 
it is clear that the index of $W_\Psi$ is exactly the same 
as the index defined in Definition \ref{theindex} and we are done
by Section 2.

\begin{remark}
It was mentioned above that the geometric index of an 
\emph{arbitrary} endomorphism must 
vanish in the case of a Cuntz algebra. This is of course obvious from the Lefschetz formula
since the $\K$-theory of Cuntz algebras vanishes rationally.  
On the other hand, it does not seem very obvious from a geometric point of view. 
This sort of thing happens in classical topology of course: one proves existence of 
fixed points by homology computations. 

\end{remark}

\end{document}